\documentclass[11pt]{amsart}
\usepackage{amsmath}
\usepackage{amscd}
\usepackage{graphicx}
\usepackage[pdftex, pdfstartview=FitH]{hyperref}
\usepackage{framed}

\newcommand{\Sys}{\mathfrak{S}}
\newcommand{\landau}{\mathcal{O}}
\newcommand{\lie}{\mathcal{L}}
\newcommand{\R}{{\mathbb R}}

\newcommand{\C}{{\mathbb C}}

\newcommand{\N}{{\mathbb N}}
\newcommand{\Z}{{\mathbb Z}}
\newcommand{\RP}{{\mathbb R}P}
\newcommand{\KP}{{\mathbb K}P}
\newcommand{\vol}{{\rm vol}}
\newcommand{\sys}{\ell_1}
\newcommand{\kernel}{\mbox{\rm Ker}\,}
\newcommand{\image}{\mbox{\rm Im}\,}
\newcommand{\dotdot}{{\, . \, . \,}}

\newcommand{\reeb}{{R_\alpha}}

\newtheorem{point}{}
\newtheorem{theorem}{Theorem}[section]

\newtheorem{lemma}[theorem]{Lemma}
\newtheorem{proposition}[theorem]{Proposition}

\theoremstyle{definition}
\newtheorem{definition}[theorem]{Definition}

\newtheorem*{question}{Question}

\newtheorem*{claim}{Claim}

\title{Contact geometry and isosystolic inequalities}
 
\author{J.C. \'Alvarez Paiva}

\address{J.C. \'Alvarez Paiva, Laboratoire Paul Painlev\'e, Bat. M2, Universit\'e des Sciences et Technologies,
59 655 Villeneuve d'Ascq, France.}

\email{juan-carlos.alvarez-paiva@math.univ-lille1.fr}

\author{F. Balacheff}

\address{F. Balacheff, Laboratoire Paul Painlev\'e, Bat. M2, Universit\'e des Sciences et Technologies,
59 655 Villeneuve d'Ascq, France.}

\email{florent.balacheff@math.univ-lille1.fr}



\begin{document}%

\begin{abstract}
A long-standing open problem in systolic geometry asks whether a Riemannian metric on the real projective space whose 
volume equals that of the canonical metric, but is not isometric to it, must necessarily carry a periodic geodesic 
of length smaller than $\pi$.  A contact-geometric reformulation of systolic geometry  and the use of canonical 
perturbation theory allow us to solve a parametric version of this problem. Namely, we show that if $g_s$ 
is a smooth volume-preserving deformation of the canonical metric and at $s=0$ the deformation is not tangent to all 
orders to trivial deformations ({\it i.e.}, to deformations of the form $\phi_s^* g_0$ for some isotopy $\phi_s$), 
then the length of the shortest periodic geodesic of the metric $g_s$ attains $\pi$ as a strict local maximum at $s=0$. 
This result still holds for complex and quaternionic projective spaces as well as for the Cayley plane. Moreover, the
same techniques can be applied to show that Zoll Finsler manifolds are the unique smooth critical points of 
the systolic volume.
\end{abstract}

\maketitle

\begin{flushright}
 \begin{small}
\parbox{4in} {{\it Pour r\'esoudre un probl\`eme nouveau, nous cherchons toujours \`a le simplifier par une s\'erie de transformations~; mais cette
simplification a un terme, car il y a dans tout probl\`eme quelque chose d'essentiel, pour ainsi dire, que toute transformation est 
impuissante \`a modifier.}} \\
--- Henri Poincar\'e
 \end{small}
\end{flushright}


\section{Introduction}

The twofold purpose of this work is to introduce contact geometry as a natural setting for the study of systolic 
inequalities---inequalities involving the shortest length of periodic geodesics and the volume of a Riemannian manifold---and to 
exploit the large symmetry group that the theory inherits through the application of canonical perturbation theory. 

The starting point of our investigations was the following concrete problem. 

\medskip
\noindent
{\bf Local systolic minimality of the round sphere.} Let $K \subset \R^3$ be a convex body with smooth boundary $\partial K$ and let 
$\sys(\partial K)$ denote the length of the shortest closed geodesics on $\partial K$ for the metric induced by the Euclidian metric. Is it true that 
$$
{\rm Area}(\partial K) \geq \frac{1}{\pi} \, \sys(\partial K)^2 
$$
whenever $K$ is sufficiently close to the unit ball ?

\medskip

To the authors's knowledge, this question appears in print for the first time in~\cite{Balacheff:2006}. There, F.~Balacheff 
shows that any infinitesimal deformation of the standard metric on the two-sphere can be realized by a smooth deformation of 
Riemannian metrics satisfying the desired inequality. Moreover, it follows from Pu's theorem~\cite{Pu:1952} that 
{\it if $K \subset \R^3$ is a centrally symmetric convex body with smooth boundary, then  
${\rm Area}(\partial K) \geq \sys(\partial K)^2/\pi$.} In this case the convex body need not be close to a ball. 

One fact that makes the non-symmetric case interesting is that the round sphere does not minimize the ratio ${\rm Area}(\partial K)/\sys(\partial K)^2 $. 
Indeed, E.~Calabi and C.~Croke (see~\cite{Croke:1988}) have remarked that if $T$ denotes the singular convex surface consisting of two identical equilateral 
triangles glued along their boundary, then
$$
\frac{{\rm Area}(T)}{\sys(T)^2} = \frac{1}{2\sqrt{3}} \, .
$$
While $T$ is not the boundary of a convex body, it can be thickened to yield convex bodies with smooth boundary for which the ratio 
${\rm Area}(\partial K)/\sys(\partial K)^2$ is strictly smaller than $1/\pi$. In fact, the singular convex surface $T$ is conjectured to realize the minimum over all Riemannian two-spheres $M$ of the ratio ${\rm Area}(M)/\sys(M)^2$. This minimum is positive by a theorem of C.~Croke \cite{Croke:1988} and its best known  lower bound is $1/32$, a result due to R.~Rotman~\cite{Rotman:2006} who improved previous lower bounds due to C.~Croke~\cite{Croke:1988}, herself in collaboration with A.~Nabutovsky~\cite{Nabutovsky-Rotman:2002}, and S.~Sabourau~\cite{Sabourau:2004}.

Another feature that adds to the interest and difficulty of the problem is the existence of {\it Zoll surfaces\/}  
(i.e., surfaces all of whose geodesics are periodic of the same minimal period). By a theorem of V.~Guillemin 
(see~\cite{Guillemin:1976}) there are plenty of non-isometric Zoll surfaces in any neighborhood of the round metric on the 
sphere. Moreover, for any such surface the ratio of its area and the square of the length of its (prime) periodic 
geodesics is equal to $1/\pi$. As a result, the round sphere has no chance of being an isolated minimum, and any 
argument based on curvature flows---none of which is known to preserve the class of Zoll surfaces---is doomed to fail.  

The local systolic minimality problem can be considered not just for the two-sphere, but for spheres and 
(real, complex, and quaternionic) projective spaces of all dimensions and, of course, for the Cayley plane. One would at least like to know whether the standard metrics on these spaces are critical points of the {\it systolic volume\/} functional 
$g \mapsto \vol(M^n,g)/\sys(M^n,g)^n$. Here $\sys(M^n,g)$ denotes the {\it systole} defined as the infimum of the lengths of periodic geodesics. We shall show that this is indeed the case not only for compact rank-one 
symmetric spaces, but for all {\it Zoll Finsler manifolds\/} (i.e., Finsler manifolds all of whose geodesics are periodic 
with the same prime period). 

\begin{definition}
A smooth one-parameter deformation $F_s$ of a closed Finsler manifold $(M,F_0)$ is said to be {\it isosystolic}  if 
the quantity $\sys(M,F_s)$ remains constant along the deformation. 
\end{definition}

\begin{theorem}\label{Zoll-Finsler}
A closed Finsler manifold $(M,F_0)$ is Zoll if and only if for every smooth isosystolic deformation $F_s$
the derivative of the function $s \mapsto \vol(M,F_s)$ vanishes at $s = 0$.
\end{theorem}

When $(M,g_0)$ is the canonical metric in a real projective space and the isosystolic deformation is required to be Riemannian, 
M.~Berger~\cite{Berger:1972} showed that $\vol(M,g_s)$ has a critical point at $s=0$. Theorem~\ref{Zoll-Finsler} improves upon Berger's 
result in that we allow deformations by (reversible and non-reversible) Finsler metrics and show that the standard metric on $\RP^n$ 
(and on any compact rank-one symmetric space) remains a critical point. Theorem~\ref{Zoll-Finsler} also implies that relatively few manifolds 
admit a {\it smooth\/} Finsler metric as a minimum of the systolic volume. Indeed, it seems to be an open question whether any manifold not 
homeomorphic to a compact rank-one symmetric space admits a Zoll metric, Riemannian or Finsler. For example, it is easy to see that the torus 
does not admit any smooth Zoll Finsler metric.

In studying the functional $F \mapsto \vol(M^n,F)/\sys(M^n,F)^n$ by perturbation techniques, we face the problem that it is not
differentiable. In Theorem~\ref{Zoll-Finsler} we bypassed this difficulty by considering smooth isosystolic deformations. However,
one of the key features of the present work is that we are able to tackle the problem head on and work with arbitrary smooth deformations,
which we sometimes normalize to be volume-preserving. Roughly speaking, our main result states that if a one-parameter deformation $F_s$ of a Zoll manifold $(M,F_0)$  does not mischievously start tangent to all orders to deformations by Zoll metrics,  the systolic volume $s \mapsto \vol(M^n,F_s)/\sys(M^n,F_s)^n$ attains a {\it strict\/} mimimum at $s=0$.

\begin{definition}
A smooth Finsler deformation $F_s$ of a Zoll manifold $(M,F_0)$ is said to be {\it formally trivial\/} if for every $m \in \N$ there exists
a deformation $F^{(m)}_s$ by Zoll Finsler metrics that has $m$-order contact with $F_s$ at $s=0$.  
\end{definition}

\begin{theorem}\label{Finsler-result}
Let $(M,F_s)$ be a smooth volume-preserving Finsler deformation of a Zoll manifold $(M,F_0)$. If the deformation is not formally trivial, then
the function $s \mapsto \sys(M,F_s)$ attains a strict local maximum at $s=0$. If, on the other hand, the deformation is 
formally trivial, then 
$$
\sys(M,F_s) = \sys(M,F_0) + \landau(|s|^k) \mbox{ for all }  k > 0 .
$$
\end{theorem}

One of the major open problems in systolic geometry is to determine whether the canonical Riemannian metric in $\RP^n$ $( n > 2)$ is a minimum of 
the systolic volume.  Specializing Theorem~\ref{Finsler-result} to the Riemannian setting and using the solution of the infinitesimal Blaschke 
conjecture by R.~Michel~\cite{Michel:1973} and C.~Tsukamoto~\cite{Tsukamoto:1981}, we obtain the following partial solution:

\begin{theorem}\label{projective-spaces}
Let $g_s$ be a smooth volume-preserving deformation of the canonical metric in the projective space $\KP^n$. If at $s=0$ 
the deformation $g_s$ is not tangent to all orders to trivial deformations ({\it i.e.}, to deformations of the form 
$\phi_s^* g_0$ for some isotopy $\phi_s$), then the length of the shortest periodic geodesic of the metric $g_s$ 
attains $\pi$ as a strict local maximum at $s=0$.
\end{theorem}

In other words, we are faced with the following non-exclusive alternatives: (1) either the length of the shortest periodic 
geodesic of the metric 
$g_s$ attains $\pi$ as a strict local maximum at $s=0$, or (2) for every $k \geq 1$ the deformation
$g_s$ is tangent to order $k$ at $s=0$ to a trivial deformation. Spheres admit non-trivial Zoll deformations so the 
situation is more delicate. However, on the two-sphere the result takes a particularly simple form.

\begin{theorem}\label{two-sphere}
Let $g_0$ be the canonical metric on the two-sphere and let $\dot{\rho} : S^2 \rightarrow \R$ be a smooth function with zero average. If $\dot{\rho}$ 
is not odd and $e^{\rho_s}g_0$ is any smooth deformation satisfying  $d\rho_s / ds \, |_{s=0} = \dot{\rho}$, then the length of the shortest periodic 
geodesic of $(S^2, e^{\rho_s}g_0)$ attains $2\pi$ as a strict local maximum at $s=0$.
\end{theorem}

This result is sharp: the main theorem of~\cite{Guillemin:1976} states that if $\dot{\rho}$ is odd, then there exists a smooth deformation 
$e^{\rho_s}g_0$ by Zoll metrics satisfying  $d\rho_s / ds \, |_{s=0} = \dot{\rho}$. The length of the shortest periodic geodesic is then constantly equal 
to $2\pi$ along the deformation.

\medskip
\noindent {\bf Plan of the paper.} Section~\ref{setup} introduces contact geometry as a natural setting for the study of 
systolic inequalities. There the reader will find the statement of our main result: a systolic-geometric characterization 
of regular contact forms. The results on Zoll Finsler manifolds stated in this introduction follow easily 
from it. Section~\ref{regular} is practically a short survey of what's known about regular contact manifolds, but it contains some new results: (1) the proof that regular contact forms are the only critical points of the systolic volume; (2) the proof that the systolic volume is a Lipschitz function in a $C^2$ 
neighborhood of any regular contact form.  In Section~\ref{Lie-transforms} we use Lie transforms and the method of averaging
to compute the normal form of deformations of regular contact forms. The proof of the main result will be found in Section~\ref{proof-main-theorem} and a theorem that includes both
Theorem~\ref{projective-spaces} and Theorem~\ref{two-sphere} as particular cases will be found in Section~\ref{applications-Riemannian-geometry}.

\medskip
\noindent {\bf Acknowledgments}
This paper has been long in the making and its results were discussed at various seminars and conferences. Along the way the 
authors have benefited  from feedback from many colleagues. We specially thank A.~Nabutovsky, I.~Babenko, K.~Cielebak, H.~Geiges,  L.~Guth, E.~Opshtein, R.~Rotman, S.~Sabourau
and F.~Schlenk. The intriguing idea that Zoll manifolds should be the critical points of {\it some\/} variational problem was passed 
on to one of us (J.-C. Alvarez Paiva) by C.~Dur\'an in a conversation long, long ago. The idea to use the perturbation 
techniques of celestial mechanics to study the geodesic flows of metrics close to the round sphere was taken from 
Poincar\'e~\cite{Poincare:1905} (\S~3) and Moser~\cite{Moser:1970}.

\section{Systolic geometry from the contact viewpoint} \label{setup}

Although there is a wonderful reference in contact geometry (\cite{Geiges:book}), for the benefit of the reader we recall the 
basic notions and in so doing adapt the presentation to our needs.

\subsection{Basic definitions}
We follow Boothby and Wang's classic paper~\cite{Boothby-Wang:1958} and define a {\it contact manifold\/} as a pair $(M, \alpha)$ consisting of a 
$(2n+1)$-dimensional manifold together with a smooth $1$-form $\alpha$ such that the top order form $\alpha \wedge d\alpha^n$ never vanishes. 

\begin{framed}
Without exception, all our contact manifolds are closed and oriented in such a 
way that $\alpha \wedge d\alpha^n > 0$.
\end{framed}

The kernel of $\alpha$ defines a field of hyperplanes in the tangent space of $M$ (a vector sub-bundle of $TM$ of co-dimension one) that is maximally 
non-integrable. Standard notation and terminology will have us denote this sub-bundle by $\xi$ and call it the {\it contact structure\/} associated 
to the contact form $\alpha$. Note that if $(M,\alpha)$ is a contact manifold and $\rho : M \rightarrow \R$ is a smooth function that never vanishes, 
the form $\rho \alpha$ is also a contact form which defines the same contact structure as $\alpha$. If we are only interested in the contact structure
associated to a contact form on $M$, we will write $(M,\xi)$. 

Diffeomorphisms of $(M,\alpha)$ that preserve the contact structure are called {\it contactomorphisms\/} or {\it contact transformations,} 
while diffeomorphisms that preserve the contact form are {\it strict contactomorphisms.} Two contact manifolds $(M_1,\alpha_1)$  and $(M_2,\alpha_2)$ 
are {\it contactomorphic\/}  if there exists a diffeomorphism $\phi : M_1 \rightarrow M_2$ and a nowhere-vanishing smooth function $\rho : M_1 \rightarrow \R$ 
such that $\phi^* \alpha_2 = \rho \alpha_1$. When $\rho \equiv 1$ (i.e.,  $\phi^* \alpha_2 = \alpha_1$), the contact manifolds are said to be 
{\it strictly contactomorphic.}

As we have defined them, contact manifolds come with a natural volume:
$$
\vol(M,\alpha) := \int_M \alpha \wedge d\alpha^n \, .
$$  
They also carry a natural vector field, the {\it Reeb vector field\/} $\reeb$, defined by the equations 
$d\alpha(\reeb,\cdot) = 0$ and $\alpha(\reeb) = 1$. The flow of the vector field $\reeb$ (remember our contact 
manifolds are all closed) is called the {\it Reeb flow\/} and its orbits are the {\it Reeb orbits.} If we are not particularly interested
in the parameterization of a Reeb orbit, we shall call it a {\it characteristic\/} (a closed characteristic if the orbit is periodic). 
Equivalently, characteristics are the $1$-dimensional leaves of the {\it characteristic distribution\/} $\kernel d\alpha$.

Our main interest is to find inequalities relating the volume of a contact manifold and the bottom of its 
{\it action spectrum:} the set of periods of its periodic Reeb orbits. 

\begin{definition}
The {\it systole\/} of a contact manifold $(M,\alpha)$, which we denote by $\sys(M,\alpha)$, is the smallest period of 
any of its periodic Reeb orbits. We define the {\it systolic volume\/} of a contact manifold $(M,\alpha)$ of 
dimension $2n+1$ as the ratio
$$
\Sys(M,\alpha) = \frac{\vol(M,\alpha)}{\sys(M,\alpha)^{n+1}} \, .
$$ 
\end{definition}

Seen as a functional on the space of contact forms inducing a given contact structure $(M,\xi)$, the systolic volume is upper 
semi-continuous. More precisely:

\begin{proposition}\label{semi-continuity}
Let $(M,\alpha)$ be a contact manifold and let $C^\infty_+(M)$ be the set of (strictly) positive smooth functions on $M$. 
The functional 
$$
\Sys :  C^\infty_+(M) \longrightarrow (0 \dotdot \infty) \\
$$
defined by $\rho \mapsto \Sys(M,\rho \alpha)$ is upper semi-continuous in the $C^1$ topology.
\end{proposition}

\proof  The proposition follows from the lower semi-continuity of the functional $\rho \mapsto \sys(M,\rho \alpha)$, which---once we pass from 
$\rho \alpha$ to its Reeb vector field at the cost of one derivative---follows in turn from a general principle: if $X$ is a nowhere-zero $C^1$
vector field on a closed manifold and all its periodic orbits have periods strictly greater than $T$, then in the space of $C^1$ vector fields
there is a $C^0$ neighborhood of $X$ where every element shares this property (see Lemma~2.3 in \cite{Palis-deMelo}). 
\qed 

The reader may have already noted that in the definition of systolic volume and in the previous proposition we are implicitly assuming the existence of 
periodic Reeb orbits on closed contact manifolds. We can bypass this thorny issue by setting $\Sys(M,\alpha) = 0$ if there are no periodic Reeb orbits, 
but the (Weinstein) conjecture is that they always exist. More importantly, their existence has been proved for all the contact manifolds that appear 
in our results. We will say more about this as we go along.

\subsection{Symmetries in contact systolic geometry}
The following simple result---whose proof is left to the reader---shows that the systolic volume is invariant under a very general class
of transformations. 

\begin{proposition}\label{invariance-systolic-volume}
If $(M_1,\alpha_1)$ and $(M_2,\alpha_2)$ are contact manifolds for which there exists a smooth map $\phi : M_1 \rightarrow M_2$ 
such that $\phi^* \alpha_2 = c\, \alpha_1 + df$, where $c$ is a non-zero constant and $f$ is a smooth function on $M_1$, then
$$
\Sys(M_1,\alpha_1) = \Sys(M_2,\alpha_2).
$$ 
\end{proposition}

From the viewpoint of perturbation theory, one advantage of working in the contact setting is that we may assume that every 
smooth deformation of a contact manifold $(M,\alpha_0)$ is of the form $\rho_s \alpha_0$, where $\rho_s$ is a smooth function 
on $M$ depending smoothly on the parameter. More precisely, {\it Gray's stability theorem\/} (see Theorem~2.2.2 
in~\cite{Geiges:book}) states that given a smooth deformation $\alpha_s$ ($s$ ranging over some compact interval), there 
exists an isotopy $\Phi_s$ such that $\Phi_s^* \alpha_s = \rho_s \alpha_0$.  In other words, we may assume that the contact 
structure stays fixed along the deformation.  

Note that if we apply a contact isotopy to a deformation of the form $\alpha_s = \rho_s \alpha_0$, the result will be another 
deformation of the same form. The measure in which we can simplify (i.e., bring to normal form) the deformation $\rho_s \alpha_0$ 
depends on our ability to construct large classes of contact isotopies. A convenient way to do this is by integrating the 
following class of vector fields: 

\begin{definition}\label{infinitesimal-ct-definition}
A vector field $X$ on a contact manifold $(M,\alpha)$ is said to be an {\it infinitesimal contact transformation\/} if 
$\lie_X \alpha = \lambda \alpha$ for some smooth function $\lambda :M \rightarrow \R$. 
\end{definition}

The following standard result (cf., Theorem~2.3.1 in~\cite{Geiges:book}) shows that infinitesimal contact transformations are 
easy to construct.

\begin{proposition}\label{contact-hamiltonian}
If $h$ is a smooth function on a contact manifold $(M,\alpha)$, there is a unique vector field $X_h$---the Hamiltonian
vector field of $h$---that satisfies the equations $\alpha(X_h)=h$ and $\lie_{X_h} \alpha = \reeb(h) \alpha$. In 
particular, $X_h$ is an infinitesimal contact transformation.  
\end{proposition}

\subsection{Main results.}\label{main-results}

\begin{definition}
A contact manifold $(M, \alpha)$ is said to be {\it regular\/} if its Reeb flow is periodic and all the Reeb orbits have 
the same prime period $\sys(M,\alpha)$.
\end{definition}

\begin{theorem}\label{main-result-1}
A contact manifold $(M,\alpha)$ is regular if and only if for every smooth isosystolic deformation the derivative of the function 
$s \mapsto \vol(M,\alpha_s)$ vanishes at $s = 0$.
\end{theorem}

\begin{definition}\label{formally-trivial-def}
A smooth deformation $\alpha_s$ of a contact form $\alpha_0$ is said to be {\it trivial\/} if there exist a smooth real-valued function 
$\lambda(s)$ and an isotopy $\Phi_s$ such that $\alpha_s = \lambda(s) \Phi_s^* \alpha_0$. A smooth deformation $\alpha_s$ is said to be 
{\it formally trivial\/} if for every $m \in \N$ there exists a trivial deformation $\alpha^{(m)}_s$ that has $m$-order contact with 
$\alpha_s$ at $s=0$.  
\end{definition}

\begin{theorem}\label{main-result-2}
Let $(M,\alpha_s)$ be a smooth deformation of a regular contact manifold $(M,\alpha_0)$. If the deformation is not formally trivial, then
the function $s \mapsto \Sys(M,\alpha_s)$ attains a strict local minimum at $s=0$. If, on the other hand, the deformation is 
formally trivial, then 
$$
\Sys(M,\alpha_s) = \Sys(M,\alpha_0) + \landau(|s|^k) \mbox{ for all }  k > 0 .
$$
\end{theorem}

The proofs of these results will take up part of Section~\ref{regular} and all of Sections~\ref{Lie-transforms} and~\ref{proof-main-theorem}. 

\smallskip
\noindent {\bf First steps in contact systolic topology.} It is natural to ask whether the {\it systolic constant\/} of a contact structure 
$(M,\xi)$, defined as
$$
\sigma(M,\xi) := \inf\{\Sys(M,\alpha) : \alpha \ \hbox{\rm is a contact form with} \ \kernel \alpha = \xi \} \, ,
$$    
is a non-trivial invariant in contact topology. Perhaps recent techniques such as embedded contact homology are sufficiently
powerful to establish the strict positivity of the systolic constant for at least some classes of contact structures. Indeed, 
a very particular case of a conjecture of M.~Hutchings (see Conjecture~8.5 in~\cite{Hutchings:2010}) implies that the systolic 
constant of the standard contact structure on the three-torus is greater than $1/2$. At this time even the simplest questions about 
the systolic constant look impossibly hard. It seems we must follow D'Alembert's advice: {\it Avancez et la foi vous viendra.}

\subsection{Applications of Theorems~\ref{main-result-1} and~\ref{main-result-2}}
We shall now show that most of the results stated in the introduction are easy consequences of the two preceding theorems. 
Before we start, let us mention a result of A.~Weinstein (see~\cite{Weinstein:1975}) that facilitates the application of 
Theorem~\ref{main-result-2}: {\it Let $\alpha_0$ be a regular contact form and let $(\alpha_s)_{s\in I}$ be a smooth deformation of $\alpha_0$ for 
$s$ ranging over some compact interval $I$.  The deformation $(\alpha_s)_{s\in I}$ is trivial if and only if the forms $\alpha_s$ are regular contact 
forms for all values of the parameter $s$.}

It is well known that geodesic flows of Riemannian and 
Finsler metrics are Reeb flows (see, for example, Theorem~1.5.2 in~\cite{Geiges:book}). The precise setup is as follows: 
through the Legendre transform, a (not necessarily reversible) Finsler metric $F$ on a manifold $N$ gives rise to a Hamiltonian 
$H$ defined in the slit cotangent bundle. The restriction $\alpha$ of the canonical one-form to the unit cotangent bundle $S^*_H N$ 
(i.e., the set of covectors where $H = 1$) is a contact form and its Reeb flow is the geodesic flow of the metric. A periodic 
Reeb orbit in $(S^*_H N, \alpha)$ projects down to a closed geodesic on $N$ whose length equals the period (and the action) of 
the orbit. In particular, $\sys(S^*_H N, \alpha)$ is the length of the shortest closed geodesic on $N$.

If the metric $F$ is Riemannian and the manifold has dimension $n$, the Riemannian volume of $(N,F)$ and the contact volume of 
$(S^*_H N, \alpha)$ are related by the equality
$$
\vol(S^*_H N, \alpha) = n! b_n \vol(N,F) \, ,
$$
where $b_n$ is the volume of the $n$-dimensional Euclidean unit ball. When the metric is Finsler, {\it we define the volume 
by the preceding equality.} This is the {\it Holmes-Thompson\/} volume of a Finsler manifold (see~\cite{Thompson:1996} 
and~\cite{Alvarez-Thompson:2004} for a detailed discussion of this definition). It is important to underline that even when 
the metric is reversible this is not the Hausdorff measure of the Finsler manifold seen as a metric space. In fact, by a 
result of C.~Dur\'an (see \cite{Duran:1998}), the Holmes-Thompson volume of a reversible Finsler manifold is strictly smaller 
than its Hausdorff measure when the metric is not Riemannian. This has the agreeable consequence that any systolic inequality 
proved for the Holmes-Thompson volume is immediately true for the Hausdorff measure, and cases of equality hold only 
for Riemannian metrics.

The factor $n! b_n$ that distinguishes the contact volume of the unit cotangent bundle from the Riemannian volume of the 
manifold may make the contact-geometric reformulation of classical isosystolic inequalities slightly unfamiliar to the initiated.
For example, Pu's isosystolic inequality~(\cite{Pu:1952}) reads as follows: {\it if $H$ is the Hamiltonian of any Riemannian 
metric on the projective plane, then $\Sys(S^*_H \RP^2, \alpha) \geq 4$.} 

The identifications of $\sys(S^*_H N, \alpha)$ as the length of the shortest closed geodesic on $(N,F)$ and of 
$\vol(S^*_H N, \alpha)$ as $n! b_n \vol(N,F)$ allow us to deduce Theorems~\ref{Zoll-Finsler} and~\ref{Finsler-result} as immediate consequences of Theorems~\ref{main-result-1} and~\ref{main-result-2}, once we remark that the Finsler metric $F$ is Zoll if and only if the restriction of the canonical one-form to the cotangent bundle $S^*_H N$ is regular. 

To end this section, we remark that our methods work better for Finsler metrics than for Riemannian metrics because 
only the former are stable under small contact perturbations. Indeed, if $(N,g)$ is a Riemannian manifold and $S^*N$ is its 
unit cotangent bundle, small $C^2$ perturbations of the contact form $\alpha$ correspond to Finsler 
perturbations of the metric that are not necessarily Riemannian (nor reversible).

\section{Regular contact manifolds}\label{regular}

Contact manifolds with periodic Reeb flows were introduced by G.~Reeb in~\cite{Reeb:1952} under the cryptic names of 
S.D.F.I. (when all the Reeb orbits have the same prime period) and S.D.F'.I. (when the prime periods are not all equal). 
Reeb's S.D.F.I ({\it syst\`emes dynamiques fibr\'es avec un invariant int\'egral}) were later studied by 
W.M.~Boothby and H.C.~Wang (see \cite{Boothby-Wang:1958} and Section~7.2 in~\cite{Geiges:book}) who coined the term 
{\it regular contact manifold}. Their main results are a structure theorem and a general construction:

\begin{theorem}[Boothby-Wang~\cite{Boothby-Wang:1958}]\label{Boothby-Wang}
A regular contact manifold $(M,\alpha)$ with systole $T := \sys(M,\alpha)$ is a principal circle bundle over a base 
manifold $B$ on which $(1/T) d\alpha$ induces an integral symplectic form $\omega$ $(\mbox{i.e., } [\omega] \in H^2(B,\Z))$. 
Moreover, $[\omega]$ is the Euler class of the circle bundle. 
\end{theorem}

\begin{theorem}\label{regular-construction}
Let $\pi : S^{2n+1} \rightarrow \C P^n$ be the Hopf fibration and let $\alpha$ be the standard contact form on the unit 
sphere $S^{2n+1} \subset \C^{n+1}$. If $B$ is a closed symplectic submanifold of $\C P^n$, the restriction of $\alpha$ to 
the submanifold $\pi^{-1}(B) \subset S^{2n+1}$ is a regular contact form. Moreover, every regular contact manifold $(M,\alpha)$ 
with $\sys(M,\alpha) = \pi$ is strictly contactomorphic to a manifold obtained by this construction.
\end{theorem}

This last result is not precisely what Boothby and Wang proved in~\cite{Boothby-Wang:1958}, but a small enhancement that 
uses the Gromov-Tischler characterization of those closed symplectic manifolds that admit a symplectic embedding in complex 
projective space (see~\cite{Gromov:1970} and~\cite{Tischler:1977}). In this form, Theorem~\ref{regular-construction} shows 
just how easy it is to construct regular contact manifolds. 

We can apply~Theorem~\ref{Boothby-Wang} to complement our characterization of the critical points of the systolic volume by 
a characterization of its critical values:

\begin{proposition}
The systolic volume of a regular contact manifold $(M,\alpha)$ is a positive integer. In particular, the systolic 
volume of a Zoll Finsler manifold is a positive integer.  
\end{proposition}

\proof
Let $\pi : M \rightarrow B$ be the Boothby-Wang fibration for $(M,\alpha)$ and assume without loss of generality that $\sys(M,\alpha) = 1$. Note that the 
fiber integration of the volume form $\alpha \wedge d\alpha^n$ yields precisely the form $\omega^n$ on $B$ and, therefore, 
$$
\vol(M,\alpha) = \int_M \alpha \wedge d\alpha^n = \int_B \omega^n = \langle [\omega^n], [B] \rangle \, .
$$ 
Since $[\omega]$ and, therefore, $[\omega^n]$ are integral cohomology classes, it follows that $\vol(M,\alpha)$ is a 
positive integer.  
\qed

When $(M,\alpha)$ is the unit sphere bundle of a Zoll Riemannian manifold $(N,g)$, this is result (and its proof) are due to 
A.~Weinstein~(\cite{Weinstein:1974}). In this case the systolic volume is what Weinstein calls $j(N,g)$, which is twice the 
{\it Weinstein integer\/} of the Zoll manifold.

The authors do not know of any criterion to determine whether a contact structure $(M,\xi)$ admits a regular contact form. 
A result of D.~Blair (see page 71 in~\cite{Blair:1976}) states that the torus $T^{2n+1}$ $(n > 0)$, provided with any 
contact structure, does not. In particular, the torus  does not admit any (smooth!) contact form that (globally or locally) 
minimizes the systolic volume.

Another natural question is whether two regular contact forms can define the same contact structure without being strictly 
contactomorphic. This is directly related to the old open problem that asks whether the geodesic flows of any two Zoll 
metrics on the same manifold are symplectically conjugate. Lastly, we mention the contact-geometric generalization of what 
C.T.~Yang called the {\it Weak Blaschke Conjecture\/} in the theory of Zoll manifolds (see \cite{Yang:1991}, \cite{Reznikov:1985}, 
and~\cite{Reznikov:1994}): 

\begin{question}
If two regular contact manifolds are contactomorphic, do they have the same systolic volume ? 
\end{question}

\subsection{Systolic non-criticality of non-regular contact manifolds}

The following result, together with Theorem~\ref{main-result-2}, proves Theorem~\ref{main-result-1}. 

\begin{theorem}\label{non-criticality}
If $(M,\alpha_0)$ is not a regular contact manifold, then there exists a smooth isosystolic deformation 
$\alpha_s$ such that the derivative of the function $s \mapsto \vol(M, \alpha_s)$ is negative at $s = 0$. Moreover, 
if $\alpha_0$ is invariant under the action of a compact group $G$, the deformation can also be chosen to be $G$-invariant.  
\end{theorem}

\proof
Let us assume that $\alpha_0$ is invariant under the action of a compact group $G$, which can be the trivial group $\{e\}$. 
Fix a $G$-invariant metric $d$ that induces the standard (manifold) topology on $M$ and for $\epsilon \geq 0$ define the set
$$
M_T(\epsilon) = \{x \in M : d(x, \varphi_T(x)) \leq \epsilon \} \ ,
$$
where $T = \sys(M,\alpha_0)$ and $\varphi_t : M \rightarrow M$ is the Reeb flow in $(M,\alpha_0)$.  Remark that the set 
$M_T(\epsilon)$ is nonempty, closed and $G$-invariant.

If $\alpha_0$ is not a regular contact form, then for all sufficiently small positive values of $\epsilon$, the set 
$M_T(\epsilon)$ is {\it properly\/} contained in $M$ and its complement is nonempty, open, and $G$-invariant. We can then 
find a smooth function $\dot{\rho} : M \rightarrow \R$ that vanishes identically on $M_T(\epsilon)$ and such that 
$$
\int_M \dot{\rho}\, \alpha \wedge d\alpha^{n} < 0 \ .
$$ 
Moreover, by averaging $\dot{\rho}$ over $G$ if necessary, we may assume that $\dot{\rho}$ is $G$-invariant.

\begin{claim}
The deformation we seek is given by $\alpha_s = (1 + s\dot{\rho}) \, \alpha_0$. 
\end{claim}
 
Since $M$ is closed, the function $1 + s\dot{\rho}$ is strictly positive for all sufficiently small values of $s$ and, 
therefore, the $\alpha_s$ are $G$-invariant contact forms for small values of $s$. Furthermore,
$$
\vol(M,\alpha_s) = \int_M (1 + s\dot{\rho})^{n+1} \, \alpha_0 \wedge d\alpha_0^n 
$$
and the derivative of $s \mapsto \vol(M,\alpha_s)$ evaluated at $s = 0$ equals
$$
(n+1) \int_M \dot{\rho}\, \alpha \wedge d\alpha^{n} < 0 \ .
$$

It remains to prove that $\alpha_s$ is an isosystolic deformation. For the rest of the proof we fix $\epsilon > 0$ such that $M_T(\epsilon)$ is 
{\it properly\/} contained in $M$ and restrict the parameter $s$ to a compact interval
$-s_0 \leq s \leq s_0$ $(s_0 > 0)$ chosen in such a way that
\begin{enumerate}
 \item the $\alpha_s$ are contact forms;
 \item if $\varphi^s_t$ denotes the Reeb flow in $(M, \alpha_s)$, then $d(\varphi^0_t(x), \varphi^s_t(x)) < \epsilon/4$ 
for all $s$ in $[-s_0 \dotdot s_0]$, for all $t$ in $[0 \dotdot T]$, and for all $x \in M$.
\end{enumerate}

The identity $\alpha_s \equiv \alpha_0$ on $M_T(\epsilon)$ has two important consequences: (1) the systoles of $(M,\alpha)$ 
are also periodic Reeb orbits of $(M,\alpha_s)$ and we have the inequality 
$$
T =  \sys(M,\alpha_0) \geq \sys(M,\alpha_s) \, ;
$$
(2) any periodic orbit of $\varphi^s$ with period strictly between $0$ and $T$ must pass through the complement of $M_T(\epsilon)$.

To complete the proof let us choose a number $b$ $(0 < b < T)$ sufficiently close to $T$ so that 
$d(x,\varphi^0_t(x)) > \epsilon/2$ for all $x$ in the complement of $M_T(\epsilon)$ and all $t \in [b \dotdot T]$. By the lower 
semi-continuity of the systole (Proposition~\ref{semi-continuity}), $\sys(M,\alpha_s) > b$ for all sufficiently small values 
of $s$. Moreover, if $s \in [-s_0 \dotdot s_0]$, then for all $x$ in the complement of $M_T(\epsilon)$ and all $t \in [b,T]$, 
we have that $d(x,\varphi^s_t(x)) > \epsilon/4$. We conclude that if the parameter $s$ is sufficiently small, the flow $\varphi^s$ has no 
periodic orbit with period strictly between $0$ and $T$.
\qed

\subsection{Lipschitz continuity of the systolic volume}

The statement of Theorem~\ref{main-result-2} takes for granted the existence of periodic Reeb orbits for every contact form
in a smooth deformation $(M,\alpha_s)$ of a regular contact manifold $(M,\alpha_0)$, at least for all sufficiently small 
values of the parameter. As A.~Banyaga remarks in~\cite{Banyaga:1990}, this follows from a result of V.~Ginzburg 
(Theorem~2 in~\cite{Ginzburg:1987}). Specializing Ginzburg's proof to the contact setting, we obtain a very simple
construction that is the key to understanding the regularity of the systolic volume functional in a neighborhood 
of a regular contact form.

\begin{theorem}[V.~Ginzburg~\cite{Ginzburg:1987}]\label{Ginzburg}
Given a regular contact manifold $(M,\alpha)$, there exists an open neighborhood  of zero $U \subset C^2(M)$ 
such that if $f \in U$, the Reeb flow of the contact form $(1 + f)\alpha$ has a periodic orbit.
\end{theorem}

\proof[Sketch of the proof] The key idea is to construct a map $A : U \rightarrow C^2(M)$ such that critical points of 
$A_f : M \rightarrow \R$ lie on periodic Reeb orbits of the contact form $(1 + f)\alpha$.

Let $\pi : M \rightarrow B$ be the Boothby-Wang fibration of the regular contact form $(M,\alpha)$. On $M$ we can easily
construct a Riemannian metric such that (1) for every $x \in M$ the Reeb vector $\reeb(x)$ is a unit vector orthogonal to the 
contact hyperplane $\xi_x$; (2) the Reeb flow of $(M,\alpha)$ acts by isometries. We call this metric
$g$ and fix it once and for all. Note that $g$ induces a submersive metric on $B$ via the projection $\pi : M \rightarrow B$
and that the contact hyperplanes are precisely the horizontal subspaces of the Riemannian submersion. 

Let $\epsilon > 0$ be smaller than the injectivity radius of $(M,g)$ and, for every $x \in M$, let $D_x = D_x(\epsilon)$ denote
the image under the exponential map of the Euclidean disc
$$
\{v_x \in \xi_x \subset T_xM : \|v_x\| < \epsilon\} .
$$ 
In a Riemannian fibration, geodesics that start horizontal stay horizontal. Therefore, if $y$ is any point in $D_x$,
there exists a unique horizontal geodesic in $D_x$ joining $y$ to $x$. We shall denote this geodesic by 
$\sigma(x,y)$ and remark that $\sigma$ depends smoothly on the choice of $x$ and $y$. 

Let $f > -1$ and denote the Reeb flow of the contact form $(1+f)\alpha$ by $\varphi_t^f$. For $f$ small enough in the 
$C^2$ norm, the system of hypersurfaces $D_x$ ($x \in M$) are Poincar\'e sections. We denote the first return time of 
$x$ to $D_x$ along the flow $\varphi_t^f$ by $\tau_f(x)$, and let
$$
\gamma_f(x) := \{\varphi_t^f(x) : 0 \leq t \leq \tau_f(x) \} .
$$
Think of $\gamma_f(x)$ as an oriented, unparameterized curve joining $x$ and $x_f := \varphi_{\tau_f(x)}(x)$ and 
complete it to a closed curve $\tilde{\gamma}_f(x)$ by adding the geodesic segment $\sigma(x,x_f)$.

The function $A_f : M \rightarrow \R$ is defined by the equation  
$$
A_f(x) = \int_{\tilde{\gamma}_f(x)} (1+f)\alpha = \int_{\gamma_f(x)} (1+f)\alpha \, ,
$$
where the last equation follows from the fact that the geodesic $\sigma(x,x_f)$ is horizontal. 

The assignment $x \mapsto \tilde{\gamma}_f(x)$ defines an embedding of $M$ into its space of unparameterized oriented 
loops, and the function $A_f$ is the composition of the action functional with this embedding.  Ginzburg shows 
in~\cite{Ginzburg:1987} that the critical points of $A_f$  correspond to the fixed points of the map $x \mapsto x_f$ and, 
therefore, to periodic Reeb orbits. 
\qed

Note that the map $f \mapsto \sys(M,(1+f)\alpha)$ can now be seen as a composition $f \mapsto A_f \mapsto \min A_f$. This is 
the key to the proof of the following result:

\begin{theorem}\label{Lipschitz}
Let $(M,\alpha)$ be regular contact manifold. There exists a small open neighborhood $U$ of zero 
in the Banach space $C^2(M)$ such that the function $f \mapsto \sys(M,(1+f)\alpha)$ is Lipschitz in $U$. In
particular, the systolic volume functional $f \mapsto \Sys(M,(1+f)\alpha)$ is Lipschitz in a small neighborhood
of zero in $C^2(M)$. 
\end{theorem}

\proof
If in the proof of Theorem~\ref{Ginzburg} we consider the Reeb vector field of $(1 + f)\alpha$ as a vector field
depending in a continuously differentiable way on the parameter $f$, we clearly see that the map $f \mapsto A_f$ is a $C^1$ 
map between Banach manifolds and, by the intermediate value theorem, it is Lipschitz in a small convex neighborhood of 
$0 \in C^2(M)$. 

Since the function $\min : C^0(M) \rightarrow \R$ is Lipschitz with constant one, it is a fortiori a Lipschitz function 
on $C^2(M)$. It follows that the composition $\sys(M,(1+f)\alpha) = \min A_f$ is Lipschitz in a small convex neighborhood of 
$0 \in C^2(M)$. 
\qed

\section{Lie transforms and the method of averaging} \label{Lie-transforms}

The origin of many papers is a simple observation on which everything hinges, and so it is with this one. 

\begin{theorem}\label{integral-of-motion}
Let $(M, \alpha)$ be a regular contact manifold. If a smooth, positive function $\rho : M \rightarrow (0 \dotdot \infty)$ 
is invariant under the Reeb flow of of $\alpha$, then $\sys(M,\rho \alpha) \leq \min \rho\, \sys(M,\alpha)$. In particular, 
$\Sys(M, \rho \alpha) \geq \Sys(M,\alpha)$ and equality holds if and only if  $\rho$ is constant.
\end{theorem}

\proof
The key observation is that if $u \in M$ is a point where this minimum is attained, then the closed characteristic of 
$(M, \alpha)$ passing through $u$ is also a closed characteristic for $(M,\rho \alpha)$. Indeed, if $\gamma(t)$ is 
any parameterization of this characteristic,
$$
d(\rho \alpha)(\dot{\gamma}(t), \cdot) = d\rho \wedge \alpha \, (\dot{\gamma}(t), \cdot) + 
\rho(\gamma(t)) d\alpha(\dot{\gamma}(t), \cdot)
$$ 
vanishes identically because $d\rho(\gamma(t))$ and $d\alpha(\dot{\gamma}(t), \cdot)$ vanish identically. 
Moreover, since $\rho(\gamma(t))$ is constantly equal to $\min \rho$, we have that the action of $\gamma$ 
in $(M,\rho \alpha)$ equals
$$
\int_\gamma \rho \alpha = \min \rho \int_\gamma \alpha = (\min \rho) \, \sys(M,\alpha) \, .
$$
It follows that $\sys(M, \rho \, \alpha) \leq (\min \rho) \, \sys(M,\alpha)$. 

To prove the second part of the theorem, note that
$$
\vol(M, \rho \, \alpha) = \int_M \rho^{n+1} \alpha \wedge d\alpha^n \geq (\min \rho)^{n+1} \vol(M,\alpha)
$$
with equality if and only if $\rho$ is constant. 
\qed

Constructing  integrals of motion (i.e., smooth functions that are invariant under the Reeb flow) in
a regular contact manifold $(M,\alpha)$ is an easy matter: if $\pi : M \rightarrow B$ is the Boothby-Wang 
fibration of $(M,\alpha)$ and $f$ is any smooth function on $B$, the function $f \circ \pi$ is an integral of 
motion. Alternatively, we can take any function $\rho$ on $M$ and construct its {\it averaged function}
$$
\bar{\rho}(x) := \frac{1}{T} \int_0^T \rho(\varphi_t(x)) dt \, ,
$$ 
where $\varphi_t$ is the Reeb flow of $(M,\alpha)$ and $T := \sys(M,\alpha)$. Clearly, $\bar{\rho}$ is an 
integral of motion.

Theorem~\ref{integral-of-motion} implies that if a deformation of a regular contact form $\alpha_0$ is of 
the form $\alpha_s = \Phi_s^* \rho_s \alpha_0$, where $\Phi_s$ is an isotopy and $\rho_s$ is a smooth one-parameter family of 
integrals of motion for the Reeb flow of $\alpha_0$, then the systolic volume of $\alpha_s$ attains a minimum at $s = 0$. 
Given the large number of isotopies and integrals of motion for a periodic Reeb flow, we could hope that every 
deformation of $\alpha_0$ is of this form and thus prove the local systolic minimality of regular contact manifolds. 
Of course, this idea does not work. However, the theory of normal forms---which we adapt to contact geometry---tells us that 
it {\it almost\/} works.

\begin{definition}
Let $\alpha_s = \rho_s \alpha_0$ be a smooth deformation of the regular contact form $\alpha_0$ and let  $k$ be a non-negative 
integer. We shall say that $\alpha_s$ is in {\it normal form to order $k$} if 
$$
\alpha_s = (1 + s\mu^{(1)} + \cdots + s^k\mu^{(k)} + s^{k+1}r_s) \alpha_0 \, 
$$ 
where the functions $\mu^{(i)}$ $(1 \leq i \leq k)$ are integrals of motion for the Reeb flow of $\alpha_0$ and the 
function $r_s$ is a smooth function on $M$ depending smoothly on the parameter $s$.
\end{definition}

\begin{theorem}\label{normal-forms}
Let $\alpha_s = \rho_s \alpha_0$ be a smooth deformation of the regular contact form $\alpha_0$. Given a non-negative 
integer $k$, there exists a contact isotopy such that $\Phi^{(k)*}_s \alpha_s$ is in normal form to order $k$.   
\end{theorem}

There are several inductive procedures to construct the contact isotopy appearing in this theorem. We shall follow
a technique known as the {\it method of Dragt and Finn\/} (see~\cite{Dragt-Finn:1976} and~\cite{Finn:1986}). The  
idea is to construct the isotopy as a composition $\phi^{(k)}_{s^k} \circ \cdots \circ \phi^{(1)}_s$, where 
$\phi^{(i)}_s$ denotes the flow of a vector field $X_i$ on $M$. 

The term {\it Lie transform\/} that is usually associated to this and other equivalent methods of computing normal forms 
has its origin in a simple remark: if $\phi_\tau$ is the flow of a vector field $X$ and $\beta$ is any smooth function 
or differential form on $M$, the Taylor expansion of $\phi^*_\tau \beta$ around $\tau = 0$ can be written in terms of 
the exponential series of the Lie derivative operator $\lie_X$:
$$
\exp(\tau\lie_X)\beta = \beta + \tau \lie_X \beta + \frac{\tau^2}{2} \lie_X(\lie_X \beta) + \cdots .
$$  
In what concerns us, this boils down to remembering that, up to terms of order $k+1$ and higher, 
$\phi_{s^k}^* \beta = \beta + s^k \lie_X \beta$.

In order to construct the vector fields $X_1$, $\dots$, $X_k$ that the method requires, we shall need the following

\begin{lemma}\label{decomposition}
If $(M,\alpha)$ is a regular contact manifold, the space of smooth functions on $M$ decomposes as a direct sum of the 
kernel and the image of the operator $f \mapsto R_{\alpha}(f)$. 
\end{lemma}

\proof
If $\varphi_t$ denotes the Reeb flow of $(M,\alpha)$ and $T := \sys(M,\alpha)$, the projection onto $\kernel \reeb$ 
associated to this  decomposition is the operator that sends a function smooth function $f$ on $M$ to its averaged 
function 
$$
\bar{f}(x) := \frac{1}{T} \int_0^T f(\varphi_t(x)) dt \, .
$$ 

If the averaged function is identically zero, it is easily checked that $f = \reeb(h)$, where
$$
h(x) = \frac{1}{T} \int_0^T t f(\varphi_t(x)) dt.
$$
Therefore, we can uniquely decompose every smooth function $f$ on $M$ as 
$\bar{f} + (f - \bar{f}) \in \kernel \reeb \oplus \image \reeb$.
\qed

\proof[Proof of Theorem~\ref{normal-forms}] 
This is a simple proof by induction. The case $k=0$ follows immediately from the Taylor expansion $\rho_s = (1 + sr_s)$ around 
$s=0$. Let us now assume that the deformation $\alpha_s$ is already in normal form to order $k-1$ and show that there is a 
contact isotopy $\Phi_s$ such that $\Phi_s^* \alpha_s$ is in contact form to order $k$.

By hypothesis $\alpha_s = (1 + s\mu^{(1)} + \cdots + s^{k-1}\mu^{(k-1)} + s^{k}r_s) \alpha_0$. Expanding $r_s$ around $s = 0$ 
we have
$$
\alpha_s = (1 + s\mu^{(1)} + \cdots + s^{k-1}\mu^{(k-1)} + s^{k}\nu + s^{k+1}r'_s) \alpha_0 .
$$ 
As in Lemma~\ref{decomposition}, write $\nu$ as $\bar{\nu} + (\nu - \bar{\nu})$  and let $h$ be a function such that
$R_{\alpha_0}(h) = -(\nu - \bar{\nu})$. 

\begin{claim}
Let $X_h$ be the Hamiltonian vector field of $h$ and let $\phi_t$ be its flow. If $\Phi_s$ is the contact isotopy defined 
by $s \mapsto \phi_{s^k}$, then $\Phi_s^* \alpha_s$ is in normal form to order $k$.  
\end{claim}

Indeed, up to terms of order $k+1$ and higher 
$$
\Phi_s^* \alpha_s = (1 +  s\mu^{(1)} + \cdots + s^{k-1}\mu^{(k-1)} + s^{k}\nu)\alpha_0  + s^k \lie_{X_h} \alpha_0. 
$$
By Proposition~\ref{contact-hamiltonian}, $\lie_{X_h} \alpha_0 = R_{\alpha_0}(h) = -(\nu - \bar{\nu})$ and, therefore,
$$
\Phi_s^* \alpha_s = (1 +  s\mu^{(1)} + \cdots + s^{k-1}\mu^{(k-1)} + s^{k}\bar{\nu} + s^{k+1}r''_s)\alpha_0 . 
$$
This gives us the desired normal form with $\mu^{(k)} = \bar{\nu}$. 
\qed

\section{Proof of the main theorem} \label{proof-main-theorem}

We recall Theorem~\ref{main-result-2} assuming, without loss of generality, that our deformations are volume-preserving. 

\begin{theorem}\label{main-result-bis}
Let $(M,\alpha_s)$ be a smooth volume-preserving deformation of a regular contact manifold $(M,\alpha_0)$. If the deformation is not formally trivial, then
the function $s \mapsto \sys(M,\alpha_s)$ attains a strict local maximum at $s=0$. If, on the other hand, the deformation is 
formally trivial, then 
$$
\sys(M,\alpha_s) = \sys(M,\alpha_0) + \landau(|s|^k) \mbox{ for all }  k > 0 .
$$ 
\end{theorem}

\proof

The key elements in the proof are~Theorem~\ref{Lipschitz} on the Lipschitz continuity of the systole in a neighborhood 
of a regular contact form, Theorem~\ref{integral-of-motion}, and the induction step in the proof of Theorem~\ref{normal-forms}. 
The proof is broken down into five steps, where the first four are technical statements followed by their proofs. 

\begin{point}
If $\alpha_s = (1 + s\nu_s + s^k r_s)\alpha_0$, where $\nu_s$ and $r_s$ are smooth functions on $M$ depending smoothly on
the parameter $s$ and $k > 1$, then 
$$
\sys(M,\alpha_s) = \sys(M,(1 + s\nu_s) \alpha_0) + \landau(|s|^k) .
$$ 
\end{point}

By Theorem~\ref{Lipschitz}, there is a $C^2$ neighborhood of zero in $C^\infty(M)$ for which the function $f \mapsto \sys(M,(1+f)\alpha_0)$ 
is Lipschitz.  It follows that for all sufficiently small values of $s$, 
$$
|\sys(M, (1 + s\nu_s + s^k r_s)\alpha_0) -  \sys(M, (1 + s\nu_s)\alpha_0)| \leq L \|s^k r_s\|_2 = \landau(|s|^k).
$$

\begin{point}
If $\mu$ is an integral of motion for the Reeb flow of $\alpha_0$ and $k$ is a positive integer, then
$$
\sys(M, (1 + s^k\mu)\alpha_0) \leq \left\{ \begin{array}{ll}
                                         (1 + s^k \min \mu) \, \sys(M,\alpha_0) & \mbox{if $s\geq 0$ or $k$ is even} \\
                                         (1 + s^k \max \mu) \, \sys(M,\alpha_0) & \mbox{if  $s < 0$ and $k$ is odd}
                                         \end{array} \right.
$$
In particular, if $(1 + s^k\mu)\alpha_0$ is a volume-preserving deformation and $\mu$ is not identically zero, then 
$\sys(M, (1 + s^k\mu)\alpha_0) < \sys(M,\alpha_0)$ for all values of $s$ different from zero.
\end{point}

By Theorem~\ref{integral-of-motion}, $\sys(M, (1 + s^k\mu)\alpha_0)$ is less than $\sys(M,\alpha_0)$ times the minimum over $M$ of the function
$(1 + s^k\mu)$. We take the minimum or the maximum of $\mu$ depending on whether $s^k$ is positive or negative. 

If the deformation is volume-preserving, the equality
\begin{eqnarray*}
\vol(M,\alpha_s) &=& \int_M (1 + s^k\mu)^{n+1} \alpha_0 \wedge d\alpha_0^n \\
                 &=&  \vol(M,\alpha_0) +(n+1)s^k \! \int_M \mu \alpha_0 \wedge d\alpha_0^n + \landau(|s|^{k+1})     
\end{eqnarray*}
implies that the integral of $\mu$ over $M$ is zero. Therefore, if $\mu$ is not identically zero, its extrema must have opposites signs. We 
conclude that the function $\sys(M, (1 + s^k\mu)\alpha_0)$ attains a strict maximum at $s=0$.

\begin{point}
If $\alpha_s = (1 + s^k\rho + s^{k+1}r_s)\alpha_0$ and $\bar{\rho}$ denotes the function obtained by averaging $\rho$ along 
the orbits of the Reeb vector field of $\alpha_0$, then there exists a contact isotopy $\Phi_s : M \rightarrow M$ such that
$$
\Phi_s^* \alpha_s = (1 + s^k\bar{\rho} + s^{k+1}r'_s) \alpha_0 \, ,
$$
where $r'_s$ is some smooth function on $M$ depending smoothly on the parameter.
\end{point}

This is precisely the induction step in Theorem~\ref{normal-forms}.

\begin{point}
Let $\alpha_s = (1 + s^k\rho + s^{k+1}r_s)\alpha_0$ be a smooth volume-preserving deformation.
\begin{enumerate}
 \item If $\bar{\rho}$ is not identically zero, then $\sys(M,\alpha_s) < \sys(M,\alpha_0)$ for all sufficiently small values of $s$ 
       different from zero.
 \item If  $\bar{\rho}$ is identically zero, then $\sys(M,\alpha_s) =\sys(M,\alpha_0) + \landau(|s|^{k+1})$.
\end{enumerate}
\end{point}

Using {\bf 3} and {\bf 1},  we have that
\begin{eqnarray*}
\sys(M,\alpha_s) &=& \sys(M,(1 + s^k\bar{\rho} + s^{k+1}r'_s) \alpha_0) \\
                 &=& \sys(M,(1 + s^k\bar{\rho}) \alpha_0) + \landau(|s|^{k+1}) .
\end{eqnarray*}
It follows immediately that if $\bar{\rho} \equiv 0$, then $\sys(M,\alpha_s) =\sys(M,\alpha_0) + \landau(|s|^{k+1})$. Assume now that $\bar{\rho}$ is not 
identically zero. By {\bf 2}, 
$$
\sys(M, (1 + s^k\bar{\rho})\alpha_0) = \left\{ \begin{array}{ll}
                                         (1 + s^k \min \bar{\rho}) \, \sys(M,\alpha_0) & \mbox{if $s\geq 0$ or $k$ is even} \\
                                         (1 + s^k \max \bar{\rho}) \, \sys(M,\alpha_0) & \mbox{if  $s < 0$ and $k$ is odd} 
                                         \end{array} \right.
$$
and, therefore, in order to prove that  $\sys(M,\alpha_s) < \sys(M,\alpha_0)$  for all sufficiently small values of $s$ different from zero, we just 
need prove that $\max \bar{\rho} > 0$ and $\min \bar{\rho} < 0$. This follows from the volume-preserving character of the deformation exactly 
as in {\bf 2.}

\begin{point}
The upshot: 
\end{point}

Let $\alpha_s = \rho_s \alpha_0$ be a smooth volume-preserving deformation of the regular contact form $\alpha_0$. 
Developing $\rho_s$ around $s=0$, we write 
$$
\alpha_s = (1 + s\rho_{(1)} + s^2 r_s) \alpha_0 \, , \mbox{\rm where } \rho_{(1)} = d\rho_s /ds \, |_{s=0} \, .
$$ 
By {\bf 4}, if the average $\bar{\rho}_{(1)}$ is not identically zero, then  $s\mapsto \sys(M,\alpha_s)$ attains a strict 
maximum at $s=0$ and we're done. If, on the other hand, $\bar{\rho}_{(1)} \equiv 0$, then, also by {\bf 4}, 
$\sys(M,\alpha_s) =\sys(M,\alpha_0) + \landau(|s|^{2})$. 

More importantly, {\bf 3} tells us that if $\bar{\rho}_{(1)} \equiv 0$, there is a contact isotopy $\Phi_s$ such that
$\Phi_s^* \alpha_s = (1 + s^2r'_s) \alpha_0$. Since this new deformation is also a smooth volume-preserving deformation of 
$\alpha_0$ and $\sys(M, \Phi_s^* \alpha_s) = \sys(M,\alpha_s)$, we can rewrite $\alpha_s =  (1 + s^2r'_s) \alpha_0$ and 
start anew. If we repeat this process an arbitrary number of times, we see that either $s \mapsto \sys(M,\alpha_s)$ attains 
a strict maximum at $s = 0$ or that $\sys(M,\alpha_s) = \sys(M,\alpha_0) + \landau(|s|^k)$ for all $k > 0$.

It also follows from the proof that if $\sys(M,\alpha_s)$ does not attain a strict maximum at $s=0$, then for any positive 
integer $k$, there exists a contact isotopy $\Phi^{(k)}_s$ and a smooth function $\nu^{(k)}_s$ on $M$ depending smoothly of 
the parameter $s$, such that $\Phi^{(k)*}_s \alpha_s = (1 + s^{k+1}\nu^{(k)}_s) \alpha_0$. In other words, either 
$\sys(M,\alpha_s)$ attains a strict maximum at $s=0$ or the deformation $\alpha_s$ is formally trivial. 
\qed

\section{Applications to Riemannian geometry} \label{applications-Riemannian-geometry}

In the results that follow we are solely concerned with Zoll Riemannian metrics 
and their Riemannian deformations.

\begin{definition}
Let $(N,g_0)$ be a Zoll manifold whose geodesic flow $\varphi_t$ has prime period $T$. A symmetric $2$-tensor $h$ on $M$ is 
said to satisfy the {\it zero-energy condition\/} (with respect to $g_0$) if the X-ray transform
$$
\hat{h}(v) = \int_0^T h(\varphi_t(v),\varphi_t(v)) \, dt
$$ 
is identically zero as a function on the unit tangent bundle of $N$.
\end{definition}

\begin{theorem}\label{Riemannian-result}
Let $(N,g_s)$ be a smooth volume-preserving deformation of a Zoll manifold $(N,g_0)$. If the deformation tensor 
$\dot{g} := dg_s/ds \, |_{s=0}$ does not satisfy the zero-energy condition, then the infimum of the  lengths of periodic geodesics of $(N,g_s)$ attains a strict local maximum at $s=0$.  
\end{theorem}

Before going into the proof of this result, we remark that Theorem~\ref{two-sphere} follows 
easily from it. 

\proof[Proof of Theorem~\ref{two-sphere}]
If the deformation of the standard metric on the two-sphere is given as $e^{\rho_s}g_0$, the zero-energy condition is 
precisely the condition that the integral of $\dot{\rho} = d\rho_s /ds \, |_{s=0}$ over any great circle be equal to zero. By 
a classic result of Funk this condition is verified if and only if $\dot{\rho}$ is an odd function on the two-sphere. We now 
apply Theorem~\ref{Riemannian-result} to conclude the proof.
\qed   

As the reader may have already guessed, the gist of the proof of Theorem~\ref{Riemannian-result} is the (obvious once you 
see it) relation between the zero-energy condition and the method of averaging in canonical perturbation theory.

\proof[Proof of Theorem~\ref{Riemannian-result}]
In order to apply the techniques of this and the previous section, we work on the unit cotangent bundle of the Zoll 
manifold $(N,g_0)$, which we denote by $S^*N$. Note that if $\alpha_0$ is the restriction of the canonical $1$-form 
to the unit cotangent bundle, $(S^*N, \alpha_0)$ is a regular contact manifold. 

For each $x \in N$, let $g_s^*(x) : T_x^*N \times T_x^*N \rightarrow \R$ be the inner product dual to $g_s(x)$. We define 
$g^*_s : T^*N \rightarrow \R$ as the function $p_x \mapsto g_s^*(x)(p_x,p_x)$ and $\rho_s$ as the restriction of 
$(g^*_s)^{-1/2}$ to $S^*N$.  The deformation of Riemannian metrics $g_s$ can be now seen as the deformation of contact 
forms $\alpha_s = \rho_s \alpha_0$. 

Note that 
$$
\dot{\rho} := d \rho_s /ds \, |_{s=0} = d (g_s^*)^{-1/2} / ds \, |_{s=0} = -(1/2) d g^*_s /ds \, |_{s=0} .
$$
Using the isomorphism $\# : T^*N \rightarrow TN$ defined by $g_0$, we see that the function $\dot{g}^* := d g^*_s /ds \, |_{s=0}$  is  
related to the deformation tensor of the deformation $g_s$ by the identity
$$
\dot{g}^*(p_x) = \dot{g}(\# p_x, \# p_x) .
$$ 
It follows that the zero-energy condition for $\dot{g}$ is satisfied if and only if the average of $\dot{\rho}$ along the 
Reeb orbits of $(S^*N,\alpha_0)$ is identically zero.  

If the deformation tensor $\dot{g}$ does not satisfy the zero-energy condition, Point~{\bf 4} in the proof of 
Theorem~\ref{main-result-bis} tells us that $\sys(M,\alpha_s)$ attains a strict local maximum at $s=0$ and this concludes 
the proof.
\qed

\medskip

We are now ready for the {\it proof of Theorem~\ref{projective-spaces}}.

\medskip

By the Michel-Tsukamoto solution of the infinitesimal Blaschke conjecture (see~\cite{Michel:1973} and~\cite{Tsukamoto:1981}), a symmetric $2$-tensor $h$ on a 
projective space satisfies the zero-energy condition with respect to the canonical metric $g_0$ if and only if there exists a vector field $X$ such that $\lie_X g_0 = h$. 
In other words, the zero energy condition characterizes the deformation tensor of (trivial) deformations of the form $\phi_s^* g_0$.

Suppose first that the smooth deformation $g_s$ of the canonical metric $g_0$ on $\KP^n$ does not agree to first order with a trivial deformation.  
The deformation tensor $\dot{g} := dg_s/ds \, |_{s=0}$ does not satisfy the zero-energy condition and the application of Theorem~\ref{Riemannian-result} concludes 
the proof in that case.

When the deformation agrees precisely to order $k-1$ $(1 < k < \infty)$ with a trivial deformation, we owe to Larry Guth the idea of reparametrizing it so that
the new deformation does not agree {\it to first order\/} with a trivial deformation.  

Composing with an isotopy if necessary, we may assume that the deformation is of the form $g_s=g_0+s^k h+\landau(s^{k+1})$, where $h$  is a symmetric two-tensor 
that does not satisfy  the zero-energy condition. The deformation
$$
\tilde{g}_s:=g_{\sqrt[k]{s}}
$$ 
is no longer smooth at $s=0$, but nevertheless the corresponding contact form on $S^\ast \KP^n$ can be written as 
$\tilde{\rho}_s\alpha_0$ with
$$
\tilde{\rho}_s=1+s\dot{\rho}+\sum_{l=1}^{k}s^{{k+l\over k}} r_{l,s},
$$
where the functions $\dot{\rho}$ and $r_{l,s}$ depend smoothly on all of their variables.
Using Points {\bf 3} and {\bf 1} in the proof of Theorem~\ref{main-result-bis},  we have that
\begin{eqnarray*}
\sys(S^* \KP^n ,\tilde{\rho}_s \alpha_0) &=& \sys(S^* \KP^n,(1 + s\dot{\rho}))+ \landau(|s|^{{k+1}\over k}) \\
                             &=& \sys(S^* \KP^n,(1 + s\bar{\dot{\rho}}))+ \landau(|s|^{{k+1}\over k}).
\end{eqnarray*}
Since $h$ does not satisfy  the zero-energy condition, the averaged function $\bar{\dot{\rho}}$ is not identically zero, 
so Point~{\bf 4} in the proof of  Theorem~\ref{main-result-bis} implies that $\sys(S^* \KP^n ,\tilde{\rho}_s \alpha_0)$ and, 
therefore, the length of the shortest periodic geodesic of the metric $g_s$ both attain a strict local maximum at $s=0$.
\qed



\end{document}